\documentstyle[12pt,amscd]{article}

\begin{document}

\newcommand{\text}[1]{\mbox{{\rm #1}}}
\newcommand{\gd}{\delta}
\newcommand{\itms}[1]{\item[[#1]]}
\newcommand{\nin}{\in\!\!\!\!\!/}
\newcommand{\sub}{\subset}
\newcommand{\cntd}{\subseteq}
\newcommand{\go}{\omega}
\newcommand{\Pa}{P_{a^\nu,1}(U)}
\newcommand{\fx}{f(x)}
\newcommand{\fy}{f(y)}
\newcommand{\gD}{\Delta}
\newcommand{\gl}{\lambda}
\newcommand{\gL}{\Lambda}
\newcommand{\half}{\frac{1}{2}}
\newcommand{\sto}[1]{#1^{(1)}}
\newcommand{\stt}[1]{#1^{(2)}}
\newcommand{\Z}{\hbox{\sf Z\kern-0.720em\hbox{ Z}}}
\newcommand{\singcolb}[2]{\left(\begin{array}{c}#1\\#2
\end{array}\right)}
\newcommand{\ga}{\alpha}
\newcommand{\gb}{\beta}
\newcommand{\gga}{\gamma}
\newcommand{\ul}{\underline}
\newcommand{\ol}{\overline}
\newcommand{\qed}{\kern 5pt\vrule height8pt width6.5pt
depth2pt}
\newcommand{\Lraro}{\Longrightarrow}
\newcommand{\Nb}{|\!\!/}
\newcommand{\NN}{{\rm I\!N}}
\newcommand{\bsl}{\backslash}
\newcommand{\gt}{\theta}
\newcommand{\op}{\oplus}
\newcommand{\C}{{\bf C}}
\newcommand{\Q}{{\bf Q}}
\newcommand{\Op}{\bigoplus}
\newcommand{\CR}{{\cal R}}
\newcommand{\tr}{\text{tr}}
\newcommand{\grr}{\omega_1}
\newcommand{\ben}{\begin{enumerate}}
\newcommand{\een}{\end{enumerate}}
\newcommand{\ndiv}{\not\mid}
\newcommand{\bab}{\bowtie}
\newcommand{\hal}{\leftharpoonup}
\newcommand{\har}{\rightharpoonup}
\newcommand{\ot}{\otimes}
\newcommand{\OT}{\bigotimes}
\newcommand{\bwe}{\bigwedge}
\newcommand{\gep}{\varepsilon}
\newcommand{\gs}{\sigma}
\newcommand{\rbraces}[1]{\left( #1 \right)}
\newcommand{\bbox}{$\;\;\rule{2mm}{2mm}$}
\newcommand{\sbraces}[1]{\left[ #1 \right]}
\newcommand{\bbraces}[1]{\left\{ #1 \right\}}
\newcommand{\OO}{_{(1)}}
\newcommand{\TT}{_{(2)}}
\newcommand{\FF}{_{(3)}}
\newcommand{\minus}{^{-1}}
\newcommand{\CV}{\cal V}
\newcommand{\CVs}{\cal{V}_s}
\newcommand{\un}{U_q(sl_n)'}
\newcommand{\on}{O_q(SL_n)'}
\newcommand{\slq}{U_q(sl_2)}
\newcommand{\olq}{O_q(SL_2)}
\newcommand{\UU}{U_{(N,\nu,\go)}}
\newcommand{\HH}{H_{n,q,N,\nu}}
\newcommand{\ct}{\centerline}
\newcommand{\bs}{\bigskip}
\newcommand{\qua}{\rm quasitriangular}
\newcommand{\ms}{\medskip}
\newcommand{\noin}{\noindent}
\newcommand{\mat}[1]{$\;{#1}\;$}
\newcommand{\raro}{\rightarrow}
\newcommand{\map}[3]{{#1}\::\:{#2}\raro{#3}}
\newcommand{\alg}{{\rm Alg}}
\def\newtheorems{\newtheorem{theorem}{Theorem}[section]
                 \newtheorem{corollary}[theorem]{Corollary}
                 \newtheorem{prop}[theorem]{Proposition}
                 \newtheorem{lemma}[theorem]{Lemma}
                 \newtheorem{definition}[theorem]{Definition}
                 \newtheorem{Proposition}[theorem]{Proposition}
                 \newtheorem{Lemma}[theorem]{Lemma}
                 \newtheorem{example}[theorem]{Example}
                 \newtheorem{Remark}[theorem]{Remark}
                 \newtheorem{claim}[theorem]{Claim}
                 \newtheorem{sublemma}[theorem]{Sublemma}
                 \newtheorem{remark}[theorem]{Remark}
                 \newtheorem{question}[theorem]{Question}
                 \newtheorem{Question}[theorem]{Question}
                 \newtheorem{conjecture}{Conjecture}[subsection]
                 \newtheorem{Definition}[theorem]{Definition}}
\newtheorems
\newcommand{\proof}{\par\noindent{\bf Proof:}\quad}
\newcommand{\dmatr}[2]{\left(\begin{array}{c}{#1}\\
                            {#2}\end{array}\right)}
\newcommand{\doubcolb}[4]{\left(\begin{array}{cc}#1&#2\\
#3&#4\end{array}\right)}
\newcommand{\qmatrl}[4]{\left(\begin{array}{ll}{#1}&{#2}\\
                            {#3}&{#4}\end{array}\right)}
\newcommand{\qmatrc}[4]{\left(\begin{array}{cc}{#1}&{#2}\\
                            {#3}&{#4}\end{array}\right)}
\newcommand{\qmatrr}[4]{\left(\begin{array}{rr}{#1}&{#2}\\
                            {#3}&{#4}\end{array}\right)}
\newcommand{\smatr}[2]{\left(\begin{array}{c}{#1}\\
                            \vdots\\{#2}\end{array}\right)}

\newcommand{\ddet}[2]{\left[\begin{array}{c}{#1}\\
                           {#2}\end{array}\right]}
\newcommand{\qdetl}[4]{\left[\begin{array}{ll}{#1}&{#2}\\
                           {#3}&{#4}\end{array}\right]}
\newcommand{\qdetc}[4]{\left[\begin{array}{cc}{#1}&{#2}\\
                           {#3}&{#4}\end{array}\right]}
\newcommand{\qdetr}[4]{\left[\begin{array}{rr}{#1}&{#2}\\
                           {#3}&{#4}\end{array}\right]}

\title{Semisimple Triangular Hopf Algebras and
Tannakian Categories}
\author{Shlomo Gelaki\\Technion-Israel
Institute of Technology\\
Department of Mathematics\\
Haifa 32000, Israel\\
email: gelaki@math.technion.ac.il}
\maketitle

\begin{abstract}
In this paper we provide a complete and
detailed proof of Theorem 2.1 from [EG1].
This theorem states that any semisimple triangular Hopf
algebra over an algebraically closed field $k$ of
characteristic $0$ is obtained from a group algebra $k[G],$ of
a unique (up to isomorphism) finite group $G,$ by twisting its
usual comultiplication in the sense of Drinfeld [Dr]. This
result is one of the key theorems in the classification of
semisimple and cosemisimple triangular Hopf algebras over
{\em any} algebraically closed field, which was obtained in
[EG4]. The proof of this theorem relies on Deligne's theorem
on Tannakian categories [De] in an essential way. We thus, in
particular, review the notions of tensor,
rigid, symmetric and Tannakian categories, discuss Deligne's
theorem, and study the category of
finite-dimensional representations of a finite-dimensional
triangular (semisimple) Hopf algebra.
\end{abstract}
\maketitle

\section{Introduction}
One of the most fundamental problems in the theory of Hopf
algebras is the classification and construction of
finite-dimensional Hopf algebras over an algebraically closed
field $k$ of characteristic $0.$ However, this problem is so
difficult that until today there exists in the literature only a
single general classification result; namely, the classification
of Hopf algebras of prime dimension. Kaplansky conjectured in 1975
that any prime dimensional Hopf algebra over $k$ is isomorphic to
the group algebra $k[\Z_p]$ [Kap], and it was only in 1994 that
this conjecture was settled by Zhu [Z]. Therefore people restrict
the problem to certain classes of finite-dimensional Hopf
algebras, e.g. to semisimple ones (these are always
finite-dimensional by a theorem of Sweedler [S]). In many senses
semisimple Hopf algebras over $k$ deserve to be considered as
"quantum" analogue of finite groups, but even so, the problem
remains extremely hard even in low dimensions. So far, proceeding
the classification of semisimple Hopf algebras over $k$ by the
dimension has not proved to be very fruitful. In fact, except for
a very few special low dimensions (see e.g. [Mon2] and references
therein), the only general classification theorems known today in
the literature are in dimension $pq$ and $p^3$ where $p,q$ are
prime numbers. In the first case it was proved that semisimple
Hopf algebras of dimension $pq$ are either group algebras or duals
of group algebras [EG7], and in the second case Masuoka proved
that there are exactly $p+1$ non-isomorphic semisimple Hopf
algebras of dimension $p^3$ which are neither group algebras nor
duals of group algebras.

A great boost to the theory of (finite-dimensional) Hopf algebras
was given by Drinfeld in the mid 80's when he invented the so
called {\em quasitriangular} Hopf algebras [Dr] for the purpose of
constructing solutions to the quantum Yang-Baxter equation that
arises in mathematical physics. Quasitriangular Hopf algebras are
the Hopf algebras whose category of finite-dimensional
representations is {\em braided rigid}, and thus they have
intriguing relationships also with low dimensional topology. A
quasitriangular Hopf algebra is called {\em triangular} if and
only if the corresponding braided rigid category is {\em
symmetric}; just like the categories of finite-dimensional
representations of a group or a Lie algebra are. In particular,
Drinfeld showed that {\em any} finite-dimensional Hopf algebra can
be embedded into a finite-dimensional quasitriangular Hopf algebra
known as its Drinfeld (or quantum) double. It is thus natural to
focus on the problem of classification and construction of
semisimple (quasi)triangular Hopf algebras. For semisimple
quasitriangular Hopf algebras the problem is still widely open.
However, the theory of semisimple triangular Hopf algebras over
$k$ is essentially closed now [EG1-EG4].

The purpose of this paper is to explain in full details how
Deligne's theorem on Tannakian categories [De] was applied in
[EG1] to prove the key structure theorem about semisimple
triangular Hopf algebras over $k.$ This theorem states that any
such Hopf algebra is obtained from a group algebra $k[G]$ of a
unique (up to isomorphism) finite group $G$ by twisting its usual
comultiplication in the sense of Drinfeld (see Theorem \ref{main}
below). This result is the key theorem in the classification of
semisimple triangular Hopf algebras over $k$ which was completed
in [EG4] (see also [G2]).

The paper is organized as follows. In Section 2 we recall the
definition of a symmetric rigid category and some of the
fundamental examples of such categories. In Section 3 we recall
the definition of a Tannakian category, and recall a Theorem of
Deligne and Milne [DM] on the structure of such categories, and a
Theorem of Deligne [De] stating an intrinsic characterization of
such categories in terms of the categorical dimensions of objects.
In Section 4 we turn to Hopf algebras and prove that the category
of finite-dimensional representations of a triangular Hopf algebra
is symmetric and rigid. In Section 5 we focus on a special class
of triangular Hopf algebras obtained from group algebras of finite
groups by twisting their usual comultiplications. In Section 6 we
explain how Deligne's theorem is applied in order to show that
these are in fact all the semisimple triangular Hopf algebras over
$k$ (see Theorem \ref{main} below). In Section 7 we conclude the
paper with brief discussions of some generalizations and
applications of Theorem \ref{main}, and some other related topics.

We refer the reader to the books [ES,Kas,Mon1,S] as references for
the general theory of Hopf algebras and quantum groups.

Throughout this paper, unless otherwise specified, the ground
field $k$ is assumed to be algebraically closed of characteristic
$0,$ and all categories are assumed to be $k-$linear and abelian.

\section{Symmetric rigid categories}
In this section we recall the definition of a symmetric rigid
category and some basic examples of such categories. The notion of
a symmetric rigid category, introduced by Mac Lane [Mac1], is an
abstract formulation of the properties of the category of
finite-dimensional vector spaces over a field.

\begin{definition}\label{tencat}
A tensor category is a sextuple $({\mathcal C}, \ot, {\bf
1},a,l,r)$ where ${\mathcal C}$ is a category, ${\bf 1}\in
\text{Ob}({\mathcal C})$ is the unit object, $\ot$ is a functor

\begin{equation}\label{ot}
\ot :{\mathcal C}\times {\mathcal C}\raro {\mathcal C},
\end{equation}

$a$ is a natural isomorphism (associativity constraint)

\begin{equation}\label{assoc}
a:\ot \circ (id_{{\mathcal C}}\times \ot)\raro \ot \circ(\ot \times
id_{{\mathcal C}})
\end{equation}
such that the diagram
\begin{equation}\label{pentagon}
\begin{CD}
W\ot (X\ot (Y\ot Z)) @ > a_{W,X,Y\ot Z}>>(W\ot X)\ot (Y\ot Z) @
>a_{W\ot X,Y,Z}>>((W\ot X)\ot Y)\ot Z\\ @ V id_{W}\ot a_{X,Y,Z}VV
& & @ A a_{W,X,Y}\ot id_{Z} AA
\end{CD}
\begin{picture}(0,0)
\thinlines \put(-3.65,-1.2){\shortstack{$(W\ot (X\ot Y))\ot Z$}}
\put(-10.9,-1.1){\vector(1,0){6.8}} \put(-8.5,-0.7){$a_{W,X\ot
Y,Z}$} \put(-15.1,-1.2){\shortstack{$W\ot ((X\ot Y)\ot Z)$}}
\end{picture}
\end{equation}

\vspace*{0.5cm}
\noin
is commutative for all $W,X,Y,Z\in
\text{Ob}({\mathcal C})$ (the Pentagon Axiom), and $l_X,r_X$ are
natural isomorphisms (unit constrains)

\begin{equation}\label{unot}
l_X:{\bf 1}\ot X\raro X,\;r_X:X\ot {\bf 1}\raro X
\end{equation}

such that the diagram

\begin{equation}\label{triangle}
\begin{CD}
X\ot ({\bf 1}\ot Y) @ > a_{X,{\bf 1},Y}>>(X\ot {\bf 1})\ot Y\\
@ V id_X\ot l_Y VV @ V r_X\ot id_Y VV\\
X\ot Y @ > id_{X\ot Y} >> X\ot Y
\end{CD}
\end{equation}

is commutative for all $X,Y\in \text{Ob}({\mathcal C}).$
\end{definition}

\begin{remark}\label{coh}
{\rm The meaning of (\ref{pentagon}) and (\ref{triangle}) is the
following. Let $X_1,\dots,X_n\in \text{Ob}({\mathcal C}).$ Then by
Mac Lane's Coherence Theorem [Mac2], any two expressions obtained
from $X_1\ot \cdots\ot X_n$ by inserting ${\bf 1}$ and brackets
are isomorphic via an isomorphism composed only from the
constrains $l,r,a$ and their inverses. }
\end{remark}

There is a natural notion of morphism between two
tensor categories.

\begin{definition}\label{tenfun}
Let $\tilde{{\mathcal C}}:=({\mathcal C}, \ot, {\bf 1},a,l,r)$ and
$\tilde{{\mathcal C}'}:=({\mathcal
C}', \ot ', {\bf 1}',a',l',r')$ be two tensor categories. A tensor
functor from
$\tilde{{\mathcal C}}$ to $\tilde{{\mathcal C}'}$ is a triple $(F,J, \varphi)$ where
$F:{\mathcal C}\raro {\mathcal C}'$ is a $k-$linear functor,
$J=\{J_{X,Y}:F(X)\ot 'F(Y)\raro
F(X\ot Y)|X,Y\in \text{Ob}({\mathcal C})\}$ is a
functorial isomorphism, and $\varphi:{\bf 1}'\raro F({\bf 1})$ is an
isomorphism,
such that the diagrams

\begin{equation}\label{tenfun1}
\begin{CD}
(F(X)\ot ' F(Y))\ot ' F(Z) @ > a'_{_{F(X),F(Y),F(Z)}}>>F(X)\ot ' (F(Y)\ot '
F(Z))\\
@ V J_{X,Y}\ot ' id_{F(Z)}VV @ V id_{F(X)}\ot ' J_{Y,Z} VV\\
F(X\ot Y)\ot ' F(Z) & & F(X)\ot ' F(Y\ot Z)\\
@ V J_{X\ot Y,Z} VV @ V J_{X,Y\ot Z} VV\\
F((X\ot Y)\ot Z) @ > F(a_{_{X,Y,Z}})>>F(X\ot (Y\ot Z))
\end{CD}
\end{equation}

\begin{equation}
\begin{CD}\label{tenfun2}
{\bf 1}'\ot ' F(X)@ > l'_{_{F(X)}}>>F(X)\\
@ V \varphi\ot ' id_{F(X)} VV @ A F(l_{_{X}}) AA\\
F({\bf 1})\ot 'F(X) @ > J_{{\bf 1},X}>>F({\bf 1}\ot X)
\end{CD}
\end{equation}
and
\begin{equation}
\begin{CD}\label{tenfun3}
F(X)\ot ' {\bf 1}' @ > r'_{_{F(X)}}>>F(X)\\
@ V id_{F(X)}\ot ' \varphi VV @ A F(r_{_{X}}) AA\\
F(X)\ot ' F({\bf 1}) @ > J_{X,{\bf 1}}>>F(X\ot {\bf 1})
\end{CD}
\end{equation}
are commutative for all $X,Y,Z\in \text{Ob}({\mathcal C}).$
\end{definition}

There is a natural notion of morphism between two
tensor functors.
\begin{definition}\label{isofun}
Let $\tilde{{\mathcal C}}:=({\mathcal C}, \ot, {\bf 1},a,l,r)$ and $\tilde{{\mathcal
C}'}:=({\mathcal C}', \ot ', {\bf 1}',a',l',r')$ be two tensor categories, and
$(F_1,J_1, \varphi_1),$ $(F_2,J_2, \varphi_2)$ two tensor functors from
$\tilde{{\mathcal C}}$ to $\tilde{{\mathcal C}'}.$ A morphism of tensor
functors
$\eta:(F_1,J_1,\varphi_1)\raro (F_2,J_2, \varphi_2)$ is a natural
transformation
$\eta:F_1\raro F_2$ such that $\eta_{{\bf 1}}\circ
\varphi_1=\varphi_2,$ and the
diagram
\begin{equation}
\begin{CD}\label{isofun1}
F_1(X)\ot 'F_1(Y) @ > J_1>>F_1(X\ot Y)\\
@ V \eta_X\ot '\eta_Y VV @ V \eta_{X\ot Y}VV\\
F_2(X)\ot 'F_2(Y) @ > J_2>>F_2(X\ot Y)
\end{CD}
\end{equation}
is commutative for all $X,Y\in \text{Ob}({\mathcal C}).$ The notions of
isomorphism and automorphism, and that of an (anti-)equivalence between
two tensor
categories are defined accordingly.
\end{definition}

The following are basic examples of tensor categories.
\begin{example}\label{tenex}
{\rm Let $k$ be any field.

1. The category $\text {Vec}_k$ of all finite-dimensional vector
spaces over $k$ is a tensor category, where $\ot =\ot _k,$ ${\bf
1}=k,$ and the constrains $a,l,r$ are the usual ones.

2. Let $G$ be an affine (pro)algebraic group. The category $\text
{Rep}_k(G)$ of all finite-dimensional algebraic representations of
$G$ over $k$ is a $k-$linear abelian tensor category. Note that
the forgetful functor $Forget:\text {Rep}_k(G)\raro \text {Vec}_k$
has a tensor structure $(Forget,J,\varphi)$ where
$J_{X,Y}=id_{X\ot Y}$ for all $X,Y\in \text{Ob}(\text
{Rep}_k(G)),$ and $\varphi =id_k.$ In fact, $G$ can be
reconstructed from $\text {Rep}_k(G)$ and $(Forget,J,\varphi).$
Indeed, $G\cong \text{Aut}^{\ot }(Forget)$ the group of tensor
automorphisms of $Forget$ (see e.g. [ES, 18.2.2]). }
\end{example}

Let $({\mathcal C}, \ot, {\bf 1},a,l,r)$ be a tensor category,
and $X\in
\text{Ob}({\mathcal
C}).$ An object $X^*\in \text{Ob}({\mathcal C})$ is said to be a {\em
(left) dual} of
$X$
if there exist two morphisms $ev_X:X^*\ot X\raro {\bf 1}$ and
$coev_X:{\bf 1}\raro
X\ot X^*$ such that
the diagrams
\begin{equation}
\begin{CD}\label{dual1}
X @ > coev_X\ot id_X>>(X\ot X^*)\ot X\\
@ A id_X AA @ A a_{X,X^*,X} AA\\
X @ < id_X\ot ev_X << X\ot (X^*\ot X)
\end{CD}
\end{equation}
and
\begin{equation}
\begin{CD}\label{dual2}
X^* @ > id_{X^*}\ot coev_X>>X^*\ot (X\ot X^*)\\
@ V id_{X^*} VV @ V a_{X^*,X,X^*} VV\\
X^* @ < ev_X\ot id_{X^*} << (X^*\ot X)\ot X^*
\end{CD}
\end{equation}
are commutative. If $X\in \text{Ob}({\mathcal C})$ has a dual
object, then it is unique up to isomorphism. If {\em any} object
of ${\mathcal C}$ has a dual object, then one can define the
(contravariant) {\em dual object functor} ${\mathcal C}\raro
{\mathcal C}$ as follows: for $X,Y\in \text{Ob}({\mathcal C})$ and
$f:X\raro Y,$ $X\mapsto X^*$ and $f\mapsto f^*$ where
$f^*:Y^*\raro X^*$ is defined by $$f^*:= l_{X^*}\circ (ev_{Y}\ot
id_{X^*})\circ (id_{Y^*}\ot f\ot id_{X^*})\circ a_{Y^*,X,X^*}\circ
(id_{Y^*}\ot coev_X)\circ (r_{Y^*})^{-1}.$$

\begin{definition}\label{rigcat}
A tensor category $({\mathcal C}, \ot, {\bf 1},a,l,r)$ is called rigid if
every object $X\in \text{Ob}({\mathcal C})$ has a dual object, and the dual object
functor is an anti-equivalence of tensor categories.
\end{definition}

Let $({\mathcal C}, \ot, {\bf 1},a,l,r)$ be a tensor category, and $\tau:{\mathcal
C}\times
{\mathcal C}\raro {\mathcal C}\times {\mathcal C}$ be the permutation functor.
\begin{definition}\label{symcat}
A symmetric category $({\mathcal C}, \ot, {\bf 1},a,l,r,c)$ is
a tensor category $({\mathcal C}, \ot, {\bf 1},a,l,r)$ equipped with a natural
isomorphism $c:\ot\raro \ot\circ\tau$ (commutativity
constraint) such that the diagram
\begin{equation}\label{hexagon}
\begin{CD}
X\ot (Y\ot Z)@ > id_X\ot c_{Y,Z}>>X\ot(Z\ot Y)\\
@ V a_{X,Y,Z}VV @ V a_{X,Z,Y} VV\\
(X\ot Y)\ot Z & & (X\ot Z)\ot Y\\
@ V c_{X\ot Y,Z} VV @ V c_{X,Z}\ot id_Y VV\\
Z\ot (X\ot Y) @ > a_{Z,X,Y}>>(Z\ot X)\ot Y
\end{CD}
\end{equation}
is commutative for all $X,Y,Z\in \text {Ob}({\mathcal C})$
(Hexagon Axiom), and
\begin{equation}\label{unitary}
c_{X,Y}\circ c_{Y,X}=id_{Y\ot X}
\end{equation}
for all $X,Y\in \text {Ob}({\mathcal C}).$
\end{definition}

There is a natural notion of morphism between two symmetric categories.
\begin{definition}\label{symfun}
Let $\tilde{{\mathcal C}}:=({\mathcal C}, \ot, {\bf 1},a,l,r,c)$ and
$\tilde{{\mathcal
C}'}:=({\mathcal C}', \ot ', {\bf 1}',a',l',r',c')$ be two symmetric categories. A
symmetric functor from $\tilde{{\mathcal C}}$ to $\tilde{{\mathcal C}'}$ is
a tensor functor $(F,J,\varphi)$ such that the diagram
\begin{equation}
\begin{CD}\label{symfun1}
F(X)\ot 'F(Y) @ > c'_{F(X)\ot ' F(Y)} >> F(Y) \ot ' F(X)\\
@ V J_{X,Y} VV @ V J_{Y,X} VV\\
F(X\ot Y) @ > F(c_{X,Y}) >> F(Y\ot X)
\end{CD}
\end{equation}
is commutative for all $X,Y\in \text{Ob}({\mathcal C}).$
\end{definition}
\begin{example}\label{symex}
{\rm Let $k$ be any field.

1. The categories $\text {Vec}_k$ and $\text {Rep}_k(G)$ are
symmetric with $c=\tau.$

2. The tensor category $\text {Supervec}_k:=\text {Rep}_k(\Z_2)$
of super-vector spaces has a symmetric structure defined by
$c_{X,Y}(x\ot y)=(-1)^{|x||y|}(y\ot x)$ where $x,y$ are
homogeneous elements and $|x|\in \{0,1\}$ is the degree of $x.$
Note that as {\em symmetric} categories $\text {Supervec}_k$ and
$\text {Rep}_k(\Z_2)$ are {\em not} equivalent.}
\end{example}

In any symmetric rigid category, there is a natural notion of dimension,
generalizing the ordinary dimension of an object in $\text {Vec}_k.$
\begin{definition}\label{dimcat}
Let $({\mathcal C}, \ot, {\bf 1},a,l,r,c)$ be a symmetric rigid category and $X\in
\text{Ob}({\mathcal C}).$ The categorical dimension $\text {dim}_{{\mathcal C}}(X)\in
\text{End}({\bf 1})$ of $X$ is the morphism
\begin{equation}\label{dimc}
\text {dim}_{{\mathcal C}}(X):{\bf 1}\stackrel{ev_X}{\longrightarrow}X\ot
X^*\stackrel{c_{X,X^*}}{\longrightarrow}X^*\ot
X\stackrel{coev_X}{\longrightarrow}{\bf 1}.
\end{equation}
\end{definition}
\begin{remark}\label{presdim}
{\rm Categorical dimensions are preserved by symmetric functors.}
\end{remark}
\section{Deligne's theorem on Tannakian categories}
Let $k$ be an algebraically closed field.
Let $\tilde{{\mathcal C}}:=({\mathcal C}, \ot, {\bf 1},a,l,r,c)$ be a $k-$linear
abelian symmetric rigid category. An exact and faithful symmetric functor
$F:\tilde{{\mathcal C}}\raro \text{Vec}_{k}$ is called a {\em fiber
functor} [DM].
\begin{definition}\label{tanak}
A $k-$linear abelian symmetric rigid category which
admits a fiber functor is called Tannakian.
\end{definition}
\begin{example}\label{tanakex}
{\rm The category $\text{Rep}_k(G)$ (see Example \ref{symex}) is
Tannakian. In this
situation the forgetful functor to $\text{Vec}_k$ is a fiber functor.
Note that
$\text{End}({\bf 1})=k.$}
\end{example}

In fact, we have the following fundamental theorem of Deligne and Milne:
\begin{theorem}[DM, Theorem 2.11]\label{delmil}
Let $\tilde{{\mathcal C}}:=({\mathcal C}, \ot, {\bf 1},a,l,r,c)$
be a Tannakian category
over $k$ such that $\text{End}({\bf 1})=k,$ and let $F$ be a fiber
functor to
$\text{Vec}_k.$ Then
$\tilde{{\mathcal C}}$ is equivalent, as a symmetric rigid category, to
$\text{Rep}_k(G)$ where $G:=\text{Aut}^{\ot }(F)$ is an
affine proalgebraic group.
\end{theorem}

The following deep theorem of Deligne gives an intrinsic characterization of
Tannakian categories.
\begin{theorem}[De, Theorem 7.1]\label{del}
Let $\tilde{{\mathcal C}}:=({\mathcal C},\ot ,{\bf 1},a,l,r,c)$ be
a $k-$linear abelian symmetric rigid category over a field $k$ of
characteristic $0,$ with $\text {End}({\bf 1})=k.$ If $\dim
_{{\mathcal C}}(X)\in \Z^+$ for all $X\in \text{Ob}({\mathcal
C}),$ then $\tilde{{\mathcal C}}$ is Tannakian.
\end{theorem}

As a corollary of Theorems \ref{delmil} and \ref{del} we have the
following special
case which we shall use later in the proof of our main theorem.
\begin{corollary}\label{cordel}
Let $k$ be an algebraically closed field of characteristic $0,$
and $({\mathcal C},\ot, {\bf 1},a,l,r,c)$ a $k-$linear abelian
symmetric rigid category with $\text {End}({\bf 1})=k,$ which is
semisimple with
finitely many irreducible objects. If categorical dimensions of objects are
non-negative integers, then there exist a finite group $G$ and an
equivalence of
symmetric rigid categories  $F:{\mathcal C}\raro \text {Rep}_k(G).$
\end{corollary}
\begin{proof} Let $G$ be the affine proalgebraic group whose existence is guaranteed
by Theorems \ref{delmil} and \ref{del}.
The semisimplicity assumption,
and the assumption that $G$ has only finite number of irreducible
representations
imply that $G$ is finite.
\end{proof}

\section{The category of finite-dimensional representations of a triangular
Hopf algebra}

In this section we relate (symmetric) tensor rigid categories and
Hopf algebras. Let $(A,m,1)$ be a finite-dimensional associative
algebra with multiplication map $m$ and unit element $1$ over a
field $k,$ and let ${\mathcal C}:=\text {Rep}_k(A)$ be the
category of finite-dimensional left $A-$modules. Clearly,
${\mathcal C}$ is a $k-$linear abelian category. In fact, the
algebra $A$ can be reconstructed from ${\mathcal C}.$ Indeed, let
$Forget$ be the forgetful functor ${\mathcal C}\raro \text
{Vec}_k,$ and let $\theta:A\raro \text {End}(Forget)$ be the map
given by $\theta(a)_V:= a_{|V},$ for all $a\in A.$ Then we have
the following:
\begin{theorem}\label{pavel}
The map $\theta$ defines an isomorphism of algebras
between $A$ and $\text {End}(Forget).$
\end{theorem}
\begin{proof} See e.g. [ES, Theorem 18.1].
\end{proof}

Suppose further that $A$ is a {\em Hopf algebra}. That
is, there exist $3$ additional structure maps
$$
\Delta:A\raro A\ot A,\;\varepsilon:A\raro
k\;\text{and}\;S:A\raro A
$$ called the
comultiplication, counit and antipode respectively,
satisfying the following
properties:
\begin{equation}\label{alghom}
\Delta, \varepsilon\;\text{are algebra homomorphisms,
and}\;S\;\text {is an algebra anti-isomorphism,}
\end{equation}
\begin{equation}\label{comult}
(\Delta\ot id_A)\circ \Delta=(id_A\ot \Delta)\circ \Delta\;\text
{(i.e.}\;\Delta\;\text {is {\em coassociative}),}
\end{equation}
\begin{equation}\label{counit}
(\varepsilon\ot id_A)\circ \Delta=(id_A\ot\varepsilon)\circ \Delta=id_A\;\text {and}
\end{equation}
\begin{equation}\label{antip}
(S\ot id_A)\circ \Delta=(id_A\ot S)\circ \Delta=\varepsilon\cdot  1.
\end{equation}

What does it say about the structure of ${\mathcal C}$?
First note that the field $k$
becomes an object
of ${\mathcal C}$ thanks to $\varepsilon$ being an algebra
homomorphism:
\begin{equation}\label{trivrep}
a\cdot x=\varepsilon(a)x\;\text{for any}\;a\in A\;\text{and}\;x\in k.
\end{equation}
Second, for
any $V,W\in \text{Ob}({\mathcal C}),$ $V\ot W$ becomes an object
of ${\mathcal C}$ thanks to $\Delta$ being an
algebra homomorphism:
\begin{equation}\label{tenrep}
a\cdot(v\ot w)=\Delta(a)\cdot(v\ot w)\;\text{for any}\;a\in A,\;v\in
V\;\text{and}\;w\in W.
\end{equation}
Third, for any $V\in \text{Ob}({\mathcal C}),$
its linear dual $V^*$ becomes an
object of ${\mathcal C}$ thanks to $S$ being an algebra
anti-homomorphism:
\begin{equation}\label{dualrepr}
(a\cdot f)(v)=f(S(a)\cdot v)\;\text{for any}\;a\in A,\;v\in V\;\text{and}
\;f\in V^*.
\end{equation}
Moreover, let $V \in \text{Ob}({\mathcal C})$ and fix dual
bases $\{v_i\}$
and $\{f_i\}$ of $V,V^*$ respectively. Then it is
straightforward to verify that
the $k-$linear maps
\begin{equation}\label{ev}
ev_{V}:V^*\ot V\raro k,\;f\ot v\mapsto f(v)
\end{equation}
and
\begin{equation}\label{coev}
coev_{V}:k\raro V\ot V^*\;\text {determined by}\;
1\mapsto\sum_i v_i \ot f_i
\end{equation}
are in fact $A-$module maps. It is now straightforward to verify
that (\ref{comult})-(\ref{coev}) say precisely that $({\mathcal
C}, \ot, k,a,l,r)$ is a tensor rigid category where $a,l,r$ are
the usual constrains (as in $\text{Vec}_k$). For instance, $a$ is
an $A-$module map thanks to (\ref{comult}), and $l,r$ are
$A-$module maps thanks to (\ref{counit}).

Suppose further that $(A,m,1,\Delta,\varepsilon,S,R)$
is {\em triangular} [Dr]. That
is, $R=\sum_i a_i\ot b_i\in A\ot A$ is invertible with
$R^{-1}=R_{21}$ ($R$ is
unitary) and the following axioms hold:
\begin{equation}\label{qt1}
(\gD \ot id_A)(R)=R_{13}R_{23},\, (id_A\ot
\gD)(R)=R_{13}R_{12}
\end{equation}
and
\begin{equation}\label{qt2}
\gD^{cop}(a)R=R\gD(a)\;for\;all\;a\in A
\end{equation}
(here $\gD^{cop}=\tau\circ \gD$ where $\tau:A\ot A\raro A\ot A$ is
the usual permutation map). Note that the two identities of
(\ref{qt1}) are equivalent since $R$ is unitary. However, when $R$
is not unitary (i.e. $(A,R)$ is {\em quasitriangular} [Dr]) they
are not equivalent. Following [R] we shall say that $(A,R)$ is
{\em minimal} when $R,$ considered as an element in $\text {Hom}_k
(A^*,A),$ defines an isomorphism.

What an extra structure does it put on ${\mathcal C}$? Well,
thanks to $R$ we can define, for any $X,Y\in \text{Ob}({\mathcal
C}),$ a $k-$linear map
\begin{equation}\label{symmrep}
c_{X,Y}:X\ot Y\raro Y\ot X,\;x\ot y\mapsto \tau(R\cdot
(x\ot y)).
\end{equation}
This is in fact an $A-$module map thanks to
(\ref{qt2}). It is now
straightforward to verify that the collection
$c:=\{c_{X,Y}|X,Y\in \text
{Ob}({\mathcal C})\}$ determines a symmetric structure on
$({\mathcal C}, \ot, k,a,l,r).$
For instance, $c_{X,Y}\circ c_{Y,X}=id_{Y,X}$ thanks
to the unitarity of $R,$
and (\ref{hexagon}) is satisfied thanks to (\ref{qt1}).

In fact, the converse holds as well. Namely, given a symmetric
rigid structure on ${\mathcal C}$ such that the forgetful functor
$F:=Forget:{\mathcal C}\raro \text{Vec}_k$ is exact, faithful and
tensor, one can construct a triangular Hopf algebra structure on
the algebra $A$ (see e.g. [ES, Theorem 18.3]). Let us sketch this
construction.

Recall from Theorem \ref{pavel} that the algebras $A$ and
$\text{End}(F)$ are isomorphic. Similarly, one can show that the
algebras $A\ot A$ and $\text{End}(F^2)$ are isomorphic, where
$F^2$ is the functor ${\mathcal C}\times {\mathcal C}\raro
\text{Vec}_k,$ defined by $(V,W)\mapsto F(V)\ot F(W).$ Now, using
the tensor structure $J$ on $F,$ define the linear map
\begin{equation}\label{tk1} \Delta:\text{End}(F)\raro
\text{End}(F^2),\; \Delta(\eta)_{V,W}:=J^{-1}_{V,W}
\circ\eta_{V\ot W}\circ J_{V,W}.
\end{equation}
Also, define the linear maps \begin{equation}\label{tk2}
\varepsilon:\text{End}(F)\raro \text{End}(F({\bf 1}))=k,\;
\varepsilon(\eta):=\eta_{F({\bf 1})}
\end{equation}
and
\begin{equation}\label{tk3}
S:\text{End}(F)\raro \text{End}(F),\;S(\eta)_V:=(\eta_{V^*})^*.
\end{equation}
Finally, set
\begin{equation}\label{tk4}
R:=\tau ((c_{A,A}(1\ot 1)).
\end{equation}
Then one can show that $(A,m,1,\Delta, \varepsilon,S,R)$ is a
triangular Hopf algebra.

Summarizing the above we state the following Tannaka-Krein type
theorem.

\begin{theorem}\label{rechop}
The assignment described above defines a bijection between:
\ben
\item symmetric rigid structures on ${\mathcal C}$ such that the
forgetful functor $Forget:{\mathcal C}\raro \text{Vec}_k$ is
exact, faithful and tensor, modulo equivalence and isomorphism of
tensor functors, and
\item triangular Hopf algebra structures on $(A,m,1),$ modulo isomorphism.
\een
\end{theorem}

 From now on we assume that $A$ is a {\em triangular} Hopf algebra.
Let us now determine the categorical dimensions in ${\mathcal C}.$
Let
\begin{equation}\label{drinelm}
u:=\sum_i S(b_i)a_i
\end{equation}
be the {\em Drinfeld element} of $(A,R).$ Drinfeld
showed [Dr] that $u$ is invertible,
\begin{equation}\label{u}
uxu^{-1}=S^2(x),\;\text {for any}\; x\in A,
\end{equation}
and
\begin{equation}\label{du}
\Delta(u)=u\ot u\;\text {(i.e.}\;u\;\text {is a {\em grouplike element})}.
\end{equation}
\begin{lemma}\label{catdim}
The categorical dimension $\text{dim}_{{\mathcal C}}(X)\in k$ of an $A-$module $X$ is
given by $\text{tr}_{|X}(u).$
\end{lemma}
\begin{proof} By (\ref{dimc}), $\text{dim}_{{\mathcal C}}(X)=\sum_i (b_i\cdot
f_i)(a_i\cdot
x_i)=\sum_i
f_i (S(b_i)a_i\cdot x_i)=\sum_i f_i (u\cdot x_i)=\text{tr}_{|X}(u)$ as desired.
\end{proof}

Suppose further that $(A,m,1,\Delta,\varepsilon,S,R)$ is {\em semisimple} (i.e. $A$
is a semisimple algebra).
\begin{lemma}\label{u1}
The Drinfeld element $u$ is central, and
\begin{equation}\label{usu}
u=S(u).
\end{equation}
\end{lemma}
\begin{proof}
By a fundamental result of Larson and Radford [LR],
$S^2=id_A,$ and hence by (\ref{u}), $u$ is central.
Now, we have $(S\ot S)(R)=R$ [Dr], so $S(u)=
\sum_i S(a_i)S^2(b_i)=\sum_i a_iS(b_i)$. This shows that
$\tr(u)=\tr(S(u))$ in
every irreducible representation of $A$. But $u$ and $S(u)$ are
central,
so they act as scalars in this representation, which proves
(\ref{usu}).
\end{proof}
\begin{lemma}\label{u2}
In particular,
\begin{equation}\label{usq}
u^2=1.
\end{equation}
\end{lemma}
\begin{proof}
By (\ref{du}), $S(u)=u^{-1},$ hence the result follows from (\ref{usu}).
\end{proof}

Let us demonstrate that it is always possible
to replace $R$ with a new $R-$matrix $\tilde R$ so that the new
Drinfeld element $\tilde u$ equals $1.$ Indeed, for any irreducible
representation $V$ of $A,$ define its parity,
$p(V)\in \Z_2$, by $(-1)^{p(V)}=u|_{V}$.
Define $\tilde R\in A\ot A$ by the condition
$\tilde R|_{V\ot W}=(-1)^{p(V)p(W)}R|_{V\ot W}.$ Set
\begin{equation}\label{ru}
R_u:=\frac{1}{2}(1\ot 1+1\ot u+u\ot 1-u\ot u).
\end{equation}
\begin{lemma}\label{u3}
$\tilde R=RR_u,$ and $(A,RR_u)$ is semisimple triangular with Drinfeld
element $1.$
\end{lemma}
\begin{proof} Straightforward.
\end{proof}

This observation allows to reduce questions about semisimple triangular
Hopf algebras over $k$ to the case when the Drinfeld element is $1.$
\begin{remark}\label{dim}
{\rm One should distinguish between the categorical dimensions of objects,
defined in any braided rigid category, and their
{\em quantum dimensions}, defined only in a ribbon category (see e.g. [Kas] for the
definition of a ribbon category).
In the diagrammatic language of [Kas] the quantum dimension
corresponds to a loop without self-crossing, and the categorical
dimension to a loop with one self-overcrossing. They may be
different numbers
for a particular irreducible object. For example, in the category of
representations of a semisimple triangular Hopf algebra $(A,R)$,
quantum dimensions (for an appropriate ribbon structure) are ordinary
dimensions,
while categorical dimensions are $u|_V\text{dim}(V)$, where
$u|_{V}$ is the scalar by which the Drinfeld element $u$ acts on $V,$
i.e. $1$ or $-1.$}
\end{remark}

\section{Triangular semisimple Hopf algebras arising
from twisting in group algebras of finite groups}

Let $(A,m,1,\gD,\varepsilon,S)$ be a Hopf algebra over $k.$ Recall
[Dr] that a {\em twist} for $A$ is an invertible element $J\in
A\ot A$ which satisfies
\begin{equation}\label{t1}
(\Delta\ot id_A)(J)(J\ot 1)=(id_A\ot \Delta)(J)(1\ot J)\;\;
and \;\;(\varepsilon\ot
id_A)(J)=(id_A\ot \varepsilon)(J)=1.
\end{equation}

Drinfeld noticed that given a twist $J$ for $A,$ one
can define a new
Hopf algebra structure
$(A^J,m,1,\Delta^J,\varepsilon,S^J)$ on the algebra $(A,m,1)$
where the coproduct
$\Delta^J$
is determined by
\begin{equation}\label{t2}
\Delta^J(x)=J^{-1}\Delta(x)J
\end{equation}
for all $x\in A,$ and the antipode $S^J$ is determined by
\begin{equation}\label{t3}
S^J(x)=Q^{-1}S(x)Q
\end{equation}
for all $x\in A,$ where $Q:=m\circ(S\ot id_A)(J).$
If $A$ is triangular with the
universal $R-$matrix $R,$ then so is $A^J,$ with
the universal $R-$matrix $R^J:=J_{21}^{-1}RJ$ (
where $J_{21}:=\tau(J)$).

If $J$ is a twist for $A$ and $x$ is an invertible
element of $A$ such that $\varepsilon(x)=1,$ then
$$J^x:=\Delta(x)J(x^{-1}\otimes x^{-1})$$
is also a twist
for $A.$ We will call the twists $J$ and $J^x$ {\em gauge
equivalent}. Observe that the
map $(A^J,R^J)\raro (A^{J^x},R^{J^x})$ determined by
$a\mapsto
xax^{-1}$ is an isomorphism of triangular Hopf
algebras. Note that if
$J$ is a twist for $A,$ then $J^{-1}$ is a twist for $A^J.$
\begin{example}\label{galg}
{\rm Let $G$ be a finite group. Then $(k[G],1\ot 1)$ is semisimple
triangular. Let $J\in k[G]\ot k[G]$ be a twist. Then
$(k[G]^J,J_{21}^{-1}J)$ is also semisimple triangular. When it is
minimal triangular we shall say that $J$ is a {\em minimal}
twist.}
\end{example}
\begin{theorem}\label{teneq}
The categories $\text{Rep}_k(A)$ and $\text{Rep}_k(A^J)$ are equivalent as
symmetric rigid categories.
\end{theorem}
\begin{proof} Let $F:\text{Rep}_k(A)\raro \text{Rep}_k(A^J)$ be the functor
defined by $F(X)=X.$ Let $\varphi=id_k$
and $J=\{J_{X,Y}\},$ where $J_{X,Y}:X\ot Y\raro X\ot Y$
is defined by $J_{X,Y}(x\ot y)=J\cdot(x\ot y).$ Then it is
straightforward to verify that
$(F,J,\varphi)$ is a symmetric functor.
\end{proof}

The main result of this section is the following uniqueness theorem:
\begin{theorem}\label{uni}
Let $G,$ $G'$ be finite groups,
$H,$ $H'$ subgroups of $G,$ $G'$ respectively, and
$J,$ $J'$ minimal twists for $k[H],$ $k[H']$ respectively.
Suppose that the triangular Hopf algebras
$(k[G]^J,J_{21}^{-1}J),$ $(k[G']^{J'}J_{21}^{'-1}J')$ are isomorphic.
Then there exists a group isomorphism $\phi:G\to G'$ such that
$\phi(H)=H'$, and $(\phi\ot \phi)(J)$ is gauge equivalent to $J'$
as twists for $k[H']$.
\end{theorem}

The rest of the section is devoted to the proof of Theorem \ref{uni}.
\begin{lemma}\label{unifunctor}
Let $\mathcal C$ be the category of finite-dimensional representations
over $k$ of a finite group, and $F_1,$ $F_2:{\mathcal C}\raro \text{Vec}_k$ be
two fiber functors. Then $F_1$ is isomorphic to $F_2.$
\end{lemma}
\begin{proof} This is a special case of [DM, Theorem 3.2]. This theorem states
that the
category of fiber functors from $\text{Rep}_k(G)$ to $\text{Vec}_k,$ $G$
an affine proalgebraic group, is
equivalent to the category of $G$-torsors over $k.$ But, $k$ is algebraically
closed, hence there exists only a unique $G$-torsor over $k.$
\end{proof}

The following corollary of this lemma answers positively Movshev's
question [Mov, Remark 1] whether any symmetric twist is trivial.
\begin{corollary}\label{symtwist} Let $G$ be a finite group, and
$J$ be a symmetric twist for $k[G]$ (i.e. $J_{21}=J$). Then $J$ is gauge
equivalent to $1\ot 1.$
\end{corollary}
\begin{proof} Let ${\mathcal C}:=\text {Rep}_k(G).$
We have two fiber functors $F_1,F_2:{\mathcal C}\raro \text{Vec}_k$ arising
from the forgetful functor; namely, the trivial one and the one defined
by $J$ respectively. By Lemma \ref{unifunctor}, $F_1,F_2$ are
isomorphic. Let $\eta:F_1\raro F_2$ be an isomorphism. By definition,
$\eta=\{\eta_{_V}:V\raro V|V\in \text{Ob}({\mathcal C})\}$ is a family of
$k-$linear
isomorphisms. By naturality, the diagram
\begin{equation}
\label{geq1}
\begin{CD}
F_1(V) @ > \eta_{_V} >>F_2(V)\\
@ V F_1(f) VV @ V F_2(f) VV\\
F_1(W) @ > \eta_{_W}>>F_2(W)
\end{CD}
\end{equation}
commutes for any two objects $V,W\in \text{Ob}({\mathcal C})$ and morphism
$f:V\raro W.$
In particular, for $V=W:=k[G],$ the left regular
representation, we get
that $f\circ \eta_{_{k[G]}}=\eta_{_{k[G]}}\circ f.$ Let
$g\in k[G]$ and
$f:=r_g:k[G]\raro k[G]$ be the right multiplication by $g.$ Then,
$r_g\circ \eta_{_{k[G]}}=\eta_{_{k[G]}}\circ r_g,$ which is equivalent
to saying that $\eta_{_{k[G]}}:k[G]\raro k[G]$ is an isomorphism of the
{\em right} regular representation. Hence, $\eta_{_{k[G]}}(a)=ax$ for
all $a\in
k[G],$ where $x:=\eta_{_{k[G]}}(1)\in k[G]$ is invertible. Moreover,
setting $V:=k[G],$ $W:=k[G]\ot k[G]$ and $f:=\Delta$ in
(\ref{geq1}) yields that
$(\Delta \circ \eta_{_{k[G]}})(1)
=(\eta_{_{k[G]\ot k[G]}}\circ \Delta)(1)$ which is equivalent to
$\Delta(x)=\eta_{_{k[G]\ot k[G]}}(1\ot 1).$

Now, by (\ref{isofun1}), the diagram
\begin{equation}
\begin{CD}\label{geq2}
F_1(k[G])\ot F_1(k[G]) @ > 1\ot 1 >>F_1(k[G]\ot k[G])\\
@ V \eta_{_{k[G]}}\ot \eta_{_{k[G]}} VV @ V \eta_{_{k[G]\ot k[G]}} VV\\
F_2(k[G])\ot F_2(k[G]) @ > J >> F_2(k[G]\ot k[G])
\end{CD}
\end{equation}
commutes. In particular, $J(\eta_{_{k[G]}}\ot \eta_{_{k[G]}})(1\ot 1)=
\eta_{_{k[G]\ot k[G]}}(1\ot 1).$ Hence, $J(x\ot x)=\Delta(x)$ as
desired.
\end{proof}
\begin{remark} {\rm Here is another proof of Corollary
\ref{symtwist} which does not use Lemma \ref{unifunctor} but uses
the results of [Mov] (see also [G2] for details). Consider the
$G$-coalgebra $A_J:=k[G]$ with coproduct $\Delta_J$ determined by
$\Delta_J(x)= (x\ot x)J,$ $x\in G,$ and its dual algebra
$(A_J)^*.$ According to [Mov], this algebra is semisimple, $G$
acts transitively on its simple ideals, and $(A_J)^*$, along with
the action of $G,$ completely determines $J$ up to gauge
transformations. Clearly, since $J$ is symmetric, this algebra is
commutative. So, it is isomorphic, as a $G$-algebra, to the
algebra of functions on a set $X$ on which $G$ acts simply
transitively. Corollary \ref{symtwist} now follows from the fact
that such a $G$-set is unique up to an isomorphism (the group $G$
itself with $G$ acting by left multiplication). Note that this
proof also uses the uniqueness of $G-$torsors (the set $X$ is a
$G-$torsor).}
\end{remark}
\begin{lemma}\label{uniqueness}
Let $G,$ $G'$ be finite groups, $J,$ $J'$ twists for $k[G],$ $k[G']$
respectively, and suppose that the triangular Hopf algebras
$(k[G]^J,J_{21}^{-1}J),$ $(k[G']^{J'},J_{21}^{'-1}J')$ are isomorphic.
Then there exists a group isomorphism $\phi:G\to G'$ such that
$(\phi\ot \phi)(J)$ is gauge equivalent to $J'.$
\end{lemma}
\begin{proof} Let $f:k[G]^J\to k[G']^{J'}$ be an isomorphism
of triangular Hopf algebras.
Then $f$ defines an isomorphism of triangular Hopf algebras
from $k[G]$ to $k[G']^{J'(f\ot f)(J)^{-1}}.$ This implies that
the element $J'(f\ot f)(J)^{-1}$ is a symmetric twist
for $k[G']$. Thus, for some invertible $x\in k[G']$ one has
$J'(f\ot f)(J)^{-1}=\Delta(x)(x^{-1}\ot x^{-1})$.
Let $\phi:=\text{Ad}(x^{-1})\circ f:k[G]\to k[G']$.
It is obvious that $\phi$ is a Hopf algebra isomorphism, so it comes
from a group isomorphism $\phi:G\to G'$.
We have $(\phi\ot\phi)(J)=\Delta(x)^{-1}J'(x\ot x)$, as desired.
\end{proof}

We can now prove Theorem \ref{uni}. By Lemma \ref{uniqueness},
it is sufficient to assume that $G'=G$, and that $J$ is gauge equivalent
to $J'$ as twists for $k[G],$ and it is enough to show that there exists
an element $a\in G$ such that $aHa^{-1}=H'$ and $(a\ot
a)J(a^{-1}\ot a^{-1})$
is gauge equivalent to $J'$ as twists for $k[H'].$

So let $x\in k[G]$ be the invertible element such that
$\Delta(x)J(x^{-1}\ot x^{-1})=J'.$
In particular, this implies that $(x\ot x)R(x^{-1}\ot x^{-1})=R',$
where $R,R'$ are the $R$-matrices corresponding to $J,J'$
respectively. By the minimality of $J,$ $J',$ we have
$xk[H]x^{-1}=k[H']$. Thus,
$$J_0:=\Delta(x)(x^{-1}\ot x^{-1})=
J'(x\ot x)J^{-1}(x^{-1}\ot x^{-1})\in k[H']^{\ot 2}.$$
It is obvious that $J_0$ is a symmetric twist for $k[H']$, so by
Corollary
\ref{symtwist},
it is gauge equivalent to $1\ot 1.$ Thus, $x=x_0a$, for some
invertible $x_0\in k[H']$,
and $a\in G$. It is clear that $aHa^{-1}=H'$, and
$\Delta(x_0^{-1})J'(x_0\ot x_0)=(a\ot a)J(a^{-1}\ot a^{-1}).$
This concludes the proof of Theorem \ref{uni}.

\section{The main theorem}

We can now state and prove our main result:
\begin{theorem}\label{main}
Let $(A,R)$ be a semisimple triangular Hopf algebra over an
algebraically closed
field $k$ of characteristic $0,$ with
Drinfeld element $u,$ and set $\tilde R:=RR_u.$ Then there exist
a finite
group $G$ and a twist $J\in k[G]\ot k[G]$ such that $(A,\tilde R)$ and
$(k[G]^J,J_{21}^{-1}J)$ are isomorphic as triangular Hopf
algebras. That is, there
exists an algebra isomorphism $\phi:k[G]\raro A$ satisfying
$\gD(\phi(a))=(\phi\ot
\phi)(J^{-1}\gD(a)J),$ $\varepsilon(\phi(a))=\varepsilon(a)$ for
all $a\in k[G],$
and $(\phi\ot \phi)(J_{21}^{-1}J)=\tilde R.$
\end{theorem}
\begin{proof} Let ${\mathcal C}$ be the category of
finite-dimensional representations over $k$ of $(A,\tilde R).$ This is
a semisimple abelian
$k-$linear category with finitely many irreducible objects,
which has a structure of a symmetric rigid category (see Section 4). In this case,
the categorical dimension $\dim_{{\mathcal C}}(V)$ of $V\in \text{Ob}({\mathcal
C})$ is equal to $\tr|_V(\tilde u).$ Since the Drinfeld element $\tilde
u$ of $(A,\tilde R)$ is $1,$ it follows that it is equal to the
ordinary dimension of $V$ as a vector space. In particular, all
categorical dimensions are non-negative integers.
In this situation we can apply Theorem \ref{cordel}.

Let $G,F$ be the finite group and functor corresponding to our
category ${\mathcal C}.$ In particular, $F$ preserves categorical
dimensions of objects, and hence their ordinary dimensions. Thus, we may
identify $V,F(V)$ as vector spaces functorially for all $V\in
\text{Ob}({\mathcal C}).$
\begin{lemma}\label{1}
There exists an algebra isomorphism $\phi:k[G]\raro A$
such that for all object $V\in \text{Ob}({\mathcal C}),$ the $G-$module
structure on $F(V)$ is given via pull back along $\phi.$
\end{lemma}
\begin{proof}
Let $\{V_i|0\le i\le m\}$ be the set of all the isomorphism
classes of irreducible representations of $A.$ Since $F$ is an
equivalence of tensor rigid categories
$\{F(V_i)=V_i|0\le i\le m\}$ is the set of all the isomorphism classes of
irreducible representations of $G.$ Since $A$ and $k[G]$ are semisimple
algebras we can fix algebra isomorphisms
$k[G]\raro \bigoplus_{i=0}^m \text{End}(V_i)$ and
$\bigoplus_{i=0}^m \text{End}(V_i)\raro A.$ This determines
an isomorphism of algebras $\phi:k[G]\raro A$ (of course, this
isomorphism is not unique). Now, by the construction of
$\phi,$ the vector space $F(V_i)=V_i$ is a $G-$module via pull back along
$\phi.$
\end{proof}

By Definition \ref{tenfun}, there exists a
family of natural $G-$module isomorphisms
$J_{V,W}:F(V)\ot F(W)\raro
F(V\ot W)$ indexed by all couples $(V,W)\in \text{Ob}({\mathcal C}\times {\mathcal
C}).$
Consider $J_{A,A}:F(A)\ot F(A)\raro F(A\ot A),$ and set
$$\tilde J:=J_{A,A}(1\ot 1)\in A\ot A.$$
\begin{lemma}\label{2}
For all $V,W\in \text{Ob}({\mathcal C})$ and $v\in
V,w\in W,$
$J_{V,W}(v\ot w)=\tilde J\cdot(v\ot w).$
\end{lemma}
\begin{proof} Consider the $A-$module maps $\bar v:A\raro V$ and
$\bar
w:A\raro W$ determined by $a\mapsto a\cdot v$ and $a\mapsto a\cdot w$
for all $a\in A$ respectively.
By naturality, the diagram
\begin{equation}
\begin{CD}\label{jact}
F(A)\ot F(A) @ > J_{A,A}>>F(A\ot A)\\
@ V F(\bar v)\ot F(\bar w) VV @ V F(\bar v\ot \bar w) VV\\
F(V)\ot F(W) @ > J_{V,W}>>F(V\ot W)
\end{CD}
\end{equation}
commutes. In particular, $(J_{V,W}\circ (F(\bar v)\ot F(\bar
w)))(1\ot 1)=(F(\bar v\ot \bar
w)\circ J_{A,A})(1\ot 1),$ which is equivalent to $J_{V,W}(v\ot
w)=J_{A,A}(1\ot 1)\cdot (v\ot w)$ as desired.
\end{proof}

In particular, the isomorphism $J_{A,A}:F(A)\ot F(A)\raro
F(A\ot A)$ is determined by $x\ot y\mapsto \tilde J(x\ot y).$
Since $1\ot 1$ is in its image, it follows that
$\tilde J$ is invertible. Hence $(\phi^{-1}\ot
\phi^{-1})(\tilde J)\in
k[G]\ot k[G]$ is invertible as well. Set
\begin{equation}\label{ct}
J:=(\phi^{-1}\ot \phi^{-1})(\tilde J)^{-1} \in k[G]\ot k[G].
\end{equation}
\begin{lemma}\label{3}
For all $a\in k[G],$ $\Delta(\phi(a))=(\phi\ot
\phi)(J^{-1}\Delta(a)J).$
\end{lemma}
\begin{proof} Since the map $J_{A,A}:F(A)\ot F(A)\raro F(A\ot
A)$ is an isomorphism
of $G-$modules, we have that $J_{A,A}(a\cdot(1\ot 1))=a\cdot J_{A,A}(1\ot
1).$ By Lemma \ref{2}, this is equivalent to
$$\tilde J\left((\phi \ot
\phi)(\Delta(a))\right)=\Delta(\phi(a))\tilde J.$$
The claim follows now after replacing $\tilde J$ by $(\phi\ot
\phi)(J^{-1}).$
\end{proof}
\begin{lemma}\label{4}
For all $a\in k[G],$ $\varepsilon(\phi(a))=
\varepsilon(a).$
\end{lemma}
\begin{proof} We first show that $(\varepsilon\ot id_A)(\tilde
J)=(id_A\ot
\varepsilon)(\tilde J)=1.$ Let $r$ denote the right unit
constraints (we use the same notation for both categories for
convenience). Then by (\ref{tenfun3}), we have that the diagram
\begin{equation}
\begin{CD}\label{uc}
F(A)\ot k @ > r_{_{F(A)}}>>F(A)\\
@ V id\ot id VV @ A F(r_{_{A}}) AA\\
F(A)\ot F(k) @ > J_{A,k}>>F(A\ot k)
\end{CD}
\end{equation}
commutes. In particular, $\left (F(r_{_{A}})\circ
J_{A,k}\circ (id\ot id)\right)(1\ot 1)=
r_{_{F(A)}}(1\ot 1)$ which implies that $(\varepsilon\ot id_A)(\tilde
J)=1.$ Similarly, $(id_A \ot \varepsilon)(\tilde J)=1.$

Now, since the map $J_{A,k}:F(A)\ot k\raro F(A\ot k)$ is an
isomorphism of $G-$modules, we have that $J_{A,k}(a\cdot(1\ot 1))=a\cdot
J_{A,k}(1\ot 1).$ By Lemma \ref{2}, this is equivalent to
$\tilde
J\cdot((\phi \ot \phi)(\Delta_0(a))\cdot (1\ot 1))=\Delta(\phi(a))\cdot
(\tilde J\cdot (1\ot 1)).$
Write $\tilde J=\sum_ix_i\ot y_i.$ Then the last equation implies that
$\sum_ix_i\phi(a_1)\varepsilon(y_i)\varepsilon(\phi(a_2))=\sum_i\phi(a)_1
x_i\varepsilon(\phi(a)_2)\varepsilon(y_i)$ which in turn (since
$\phi$ is an isomorphism) implies that $\sum
a_1\varepsilon(\phi(a_2))=a,$ and the result follows.
\end{proof}
\begin{lemma}\label{5}
$J$ is a twist for $k[G].$
\end{lemma}
\begin{proof} Let $a$ denote the associativity constraints in the
categories ${\mathcal C}$ and $\text{Rep}_k(G)$ (we use the same
notation for both
categories for convenience). By (\ref{assoc}), the diagram
\begin{equation}
\begin{CD}\label{assoc2}
(F(A)\ot F(A))\ot F(A) @ > a_{_{F(A),F(A),F(A)}}>>F(A)\ot (F(A)\ot
F(A))\\
@ V J_{A,A}\ot id_{F(A)}VV @ V id_{F(A)}\ot J_{A,A} VV\\
F(A\ot A)\ot F(A) & & F(A)\ot F(A\ot A)\\
@ V J_{A\ot A,A} VV @ V J_{A,A\ot A} VV\\
F((A\ot A)\ot A) @ > F(a_{_{A,A,A}})>>F(A\ot (A\ot A))
\end{CD}
\end{equation}
commutes. In particular,
\begin{eqnarray*}
\lefteqn {\left(F(a_{_{A,A,A}})\circ J_{A\ot A,A}\circ
(J_{A,A}\ot id_{F(A)})\right)(1\ot 1\ot 1)}\\
&=& \left(J_{A,A\ot A}\circ (id_{F(A)}\ot J_{A,A})\circ
a_{_{F(A),F(A),F(A)}}\right)(1\ot 1\ot 1),
\end{eqnarray*}
which is equivalent to
$(\Delta\ot id_A)(\tilde J)\tilde J_{12}=(id_A\ot \Delta)(\tilde
J)\tilde J_{23}.$ Write $J^{-1}=\sum_ix_i\ot y_i.$ Substitute
$(\phi\ot \phi)(J)^{-1}$ for $\tilde J$ in the last equation, and use
Lemma \ref{3} to get
$$(\phi\ot \phi\ot \phi)\left((J^{-1}\Delta_0(x_i)J\ot
y_i)(J^{-1}\ot 1)\right)=
(\phi\ot \phi\ot \phi)\left((x_i\ot
J^{-1}\Delta_0(y_i)J)(1\ot J^{-1})\right).$$
Since $\phi$ is an isomorphism, this is equivalent to saying that $J$
satisfies the first part of (\ref{t1}).

Now, we already showed in the proof of Lemma \ref{4} that
$$(\varepsilon\ot id_A)(\tilde J)=(id_A\ot \varepsilon)(\tilde J)=1.$$
Thus, the second part
of (\ref{t1}) follows from Lemma \ref{3} after replacing
$\tilde J$ with
\linebreak $(\phi\ot \phi)(J^{-1}).$
\end{proof}

By Lemmas \ref{1}-\ref{5}, $(k[G]^J,J_{21}^{-1}J)$ is a triangular
semisimple Hopf algebra, and the map $\phi:k[G]^J\raro A$ is an
isomorphism of Hopf algebras. Finally, let $c$ denote the
commutativity constrains in the categories ${\mathcal C}$ and
$\text{Rep}_k(G)$ (again, we use the same notation for both
categories for convenience). Recall that $F$ is in particular a
symmetric functor.
\begin{lemma}\label{6}
$(\phi\ot \phi)(J_{21}^{-1}J)=\tilde R.$
\end{lemma}
\begin{proof} By Definition \ref{symfun}, the diagram
\begin{equation}
\begin{CD}\label{braid}
F(A)\ot F(A) @ > J_{A,A}>>F(A\ot A)\\
@ V c_{_{F(A),F(A)}} VV @ V F(c_{_{A,A}}) VV\\
F(A)\ot F(A) @ > J_{A,A}>>F(A\ot A)
\end{CD}
\end{equation}
commutes. In particular, $(F(c_{_{A,A}})\circ
J_{A,A})(1\ot 1)=(J_{A,A}\circ c_{_{F(A),F(A)}})(1\ot 1).$ Therefore,
$\tilde J=\tilde R_{21}\tilde J_{21},$ which is equivalent to the desired
result.
\end{proof}

This completes the proof of the theorem.
\end{proof}

\begin{remark}\label{sv} {\rm As seen from Remark \ref{dim}, if
$u\ne 1$, then the
category of representations of $(A,R)$ is equivalent to the category
of representations of some group as a rigid tensor category
but {\em not} as a symmetric category. This was the reason for passing
from $R$ to $\tilde R$. It is easy to see that as a symmetric
rigid category, the category of representations of $(A,R)$ is
equivalent to the category of representations of $G$ on
super-vector spaces,
such that $u$ acts by $1$ on the even part and as $-1$ on the
odd part.
For example, if $A:=k\Z_2$ with $u$ as the generator of $\Z_2,$ and $R:=R_u,$
then the category of representations is just the
category of super-vector spaces.}
\end{remark}
\begin{corollary}\label{del2}
Let $(A,R)$ be a semisimple triangular Hopf algebra over $k$
with Drinfeld element $1.$ Then there exist finite groups
$H\subseteq G,$ and a minimal twist $J\in k[H]\ot k[H]$ such
that $(A,R)\cong (k[G]^J,J_{21}^{-1}J)$ as triangular Hopf
algebras. Furthermore, the data $(G,H,J)$ is unique up to
isomorphism of groups and gauge equivalence of twists.
\end{corollary}
\begin{proof} Let $(A_R,R)$ be the minimal triangular sub Hopf algebra of
$(A,R).$ By Theorem \ref{main}, there exist a finite group $H$ and
a minimal twist $J$ for $k[H]$ such that $(A_R,R)\cong
(k[H]^J,J_{21}^{-1}J)$ as triangular Hopf algebras. We may as well
assume that $(A_R,R)=(k[H]^J,J_{21}^{-1}J).$ Let
$f:(k[H]^J,J_{21}^{-1}J)\raro (A,R)$ be the inclusion map. Then
$f:(k[H],1\ot 1)\raro (A^{J^{-1}},J_{21}RJ^{-1})$ is an injective
morphism of triangular Hopf algebras as well. In particular,
$J_{21}RJ^{-1}=1\ot 1$ which is equivalent to $R=J_{21}^{-1}J.$
Moreover, since $(A^{J^{-1}},1\ot 1)$ is triangular, $A^{J^{-1}}$
is cocommutative. Therefore there exists a finite group $G$ such
that $A^{J^{-1}}=k[G].$ Hence, $(A,R)=(k[G]^J,J_{21}^{-1}J).$
Since $k[H]$ is a sub Hopf algebra of $k[G],$ $H$ is a subgroup of
$G.$ Finally, the uniqueness follows from Theorem \ref{uni}.
\end{proof}

\section{Concluding remarks}
We conclude the paper with some related remarks and references. A
full and detailed discussion of remarks 1-5 can be found in [G2].
\begin{remark}\label{conrem}
{\rm

1. The analogue of Theorem \ref{main} in positive characteristic was
obtained in
[EG4]. In this
situation one has to further assume that $A$ is also cosemisimple (i.e.
its dual
Hopf algebra
$A^*$ is semisimple). The result that any semisimple cosemisimple
triangular Hopf
algebra over any algebraically closed field $k$ is obtained from $k[G]$
for some finite group $G,$ by twisting its usual comultiplication, follows then
using the lifting functor from positive characteristic to
characteristic
$0$ [EG5], the classification of twists for group algebras in characteristic $0$
[EG4], and Theorem \ref{main}.

2. In fact much more can be said on the structure of semisimple
cosemisimple triangular Hopf algebras over any algebraically
closed field $k.$ A complete classification of these objects was
obtained, using Theorem \ref{main} and the theory of Movshev
[Mov], in [EG4]. It was proved there that the set of isomorphism
classes of semisimple cosemisimple triangular Hopf algebras of
dimension $N$ over $k$ is in bijection with the set of isomorphism
classes of quadruples $(G,H,V,u)$ where $G$ is a finite group of
order $N,$ $H$ is a subgroup of $G,$ $V$ is an irreducible
projective representation of $H$ of dimension $|H|^{1/2}$ and
$u\in Z(G)$ is of order $\le 2.$ In particular, it follows that
$H$ is a quotient of a central type group and hence solvable by a
result of Howlett and Isaacs [HI] (it is interesting to note that
this result was proved using the classification of finite simple
groups). In [EG2] we proved that any finite group with a bijective
$1-$cocycle to an abelian group gives rise to such a group, and
hence to a minimal triangular semisimple Hopf algebra. As a
corollary to the classification theorem we were able to prove that
any semisimple cosemisimple triangular Hopf algebra has a
non-trivial grouplike element [EG4]. In [EGGS] it was shown that
this is false for semisimple quasitriangular Hopf algebras. A
counterexample was constructed there using Mathieu simple group
$M_{24}.$

3. A famous conjecture of Kaplansky from 1975 [Kap] states that
the dimension of an irreducible representation of a semisimple
Hopf algebra over $k$ divides the dimension of $A.$ This
conjecture was proved, when $A$ is quasitriangular, in [EG1] using
the theory of modular categories. In the triangular case the
conjecture also follows from Theorem \ref{main} by a famous
theorem of Frobenius. However, the conjecture does not follow
easily for the duals of semisimple triangular Hopf algebras
(semisimple {\em cotriangular} Hopf algebras). In [EG3], we used
Theorem \ref{main} to describe the representation theory of
semisimple cotriangular Hopf algebras, and in particular proved
that the conjecture holds for them as well.

4. There exist finite-dimensional triangular Hopf algebras which
are {\em not} semisimple, so in particular are not obtained by
twisting of group algebras. The simplest example is Sweedler's
$4-$dimensional Hopf algebra $A$ [S]. It is generated as an
algebra by $a,x,$ subject to the relations $a^2=1,$ $x^2=0$ and
$ax=-xa.$ Its coalgebra structure is determined by $\gD(a)=a\ot a$
and $\gD(x)=x\ot 1+a\ot x.$ This Hopf algebra is {\em pointed} but
not semisimple since it has only two non-isomorphic irreducible
representations, both of dimension $1.$ However, $A$ admits
infinitely many minimal triangular structures and one nonminimal
triangular structure (i.e. $R_a$). See [G1] for a study of pointed
non-semisimple finite-dimensional triangular Hopf algebras. In
particular, this paper gives a classification of minimal
triangular pointed Hopf algebras. The crucial property of such
Hopf algebras is that they are also {\em basic}; i.e. all their
irreducible representations are $1-$dimensional. This implies that
they share the {\em Chevalley property} with semisimple Hopf
algebras; namely, their radical is a Hopf ideal and so their
semisimple part is a Hopf algebra itself (see [AEG]). In
representation-theoretic terms this means that the tensor product
of any two irreducible representations is completely reducible.
The fact that all known examples of finite-dimensional triangular
Hopf algebras have the Chevalley property was the motivation to
study such Hopf algebras in [AEG] and [EG8]. In [AEG] Theorem
\ref{main} was used to prove that a finite-dimensional triangular
Hopf algebra has the Chevalley property if and only if it is
obtained from a super-group algebra of a finite super-group after
twisting its comultiplication. In [EG8] the main result of [EG4]
(see 2) was generalized and it was proved that the isomorphism
classes of finite-dimensional triangular Hopf algebras with the
Chevalley property are in bijection with the isomorphism classes
of septuples $(G,W,H,Y,B,V,u)$ where $G$ is a finite group, $W$ is
a finite-dimensional representation of $G,$ $H$ is a subgroup of
$G,$ $Y$ is an $H-$invariant subspace of $W,$ $B$ is an
$H-$invariant nondegenerate element in $S^2Y,$ $V$ is an
irreducible projective representation of $H$ of dimension
$|H|^{1/2}$ and $u\in Z(G)$ is of order $\le 2$ acting by $-1$ on
$W.$

5. A natural question arising from Theorem \ref{main} is what
happens in the infinite-dimensional case. Suppose $(A,R)$ is an
infinite-dimensional cotriangular Hopf algebra over $k.$ This is
equivalent to saying that the category $\text{Corep}_k(A)$ of all
the corepresentations of $A$ is symmetric rigid. We can thus ask
for a characterization of such Hopf algebras which are obtained by
twisting (this time the usual {\em multiplication} of) the
function algebra on some affine proalgebraic group. An answer to
this question is given in [EG6]. It turns out that a necessary and
sufficient condition for this to be the case is that $A$ is
pseudoinvolutive (i.e. $\tr(S^2|_{C})=\dim(C)$ for any finite
dimensional subcoalgebra of $A$). In the finite-dimensional case,
pseudoinvolutivity is equivalent to involutivity ($S^2=id_A$) ,
which in turn is equivalent to semisimplicity by a result of
Larson and Radford. Thus this result generalizes Theorem
\ref{main}. The proof also uses Deligne's theorem (Theorem
\ref{del}). The hard direction is to prove that pseudoinvolutivity
implies that the symmetry of $\text{Corep}_k(A)$ can be modified
so that the categorical dimensions of objects will become
non-negative integers.

6. As we have indicated above, if two finite groups have the same
symmetric tensor categories of representations over $\C,$ then
they are isomorphic. This raises the following natural question:
When do two finite groups $G_1,G_2$ have the same tensor
categories of representations over $\C$ (without regard for the
commutativity constraint)? Two groups with such property are
called {\em isocategorical} [EG9]. In [EG9] the theory of
triangular semisimple Hopf algebras was used to completely answer
this question. Namely, a classification of groups isocategorical
to a given group $G$ was given. In particular, it was showed that
if $G$ has no nontrivial normal abelian subgroups of order
$2^{2m}$ then any group isocategorical to $G$ must actually be
isomorphic to $G.$ On the other hand, an example of two groups
which are isocategorical but not isomorphic was given: The affine
symplectic group of a vector space over the field of two elements,
and an appropriate "affine pseudosymplectic group" introduced by
R. Griess (containing the "pseudosymplectic group" of A. Weil).
Also the notion of isocategorical groups was applied to answer the
question: When are two triangular semisimple Hopf algebras
isomorphic as Hopf algebras?}
\end{remark}

\noin
{\bf Acknowledgment} The work described in this paper
is all joint with Pavel Etingof from MIT, whom I am grateful to
for his help in reading the manuscript. I am also grateful to Bar
Ilan University for its support.


\begin{thebibliography}{[EGGS] }
\bibitem
[AEG]{aeg} N. Andruskiewitsch, P. Etingof and S. Gelaki,
Triangular Hopf Algebras with the Chevalley Property, {\em
Michigan Journal of Mathematics}, to appear, math.QA/0008232.
\bibitem
[De]{de} P. Deligne, Categories Tannakiennes, In The Grothendick
Festschrift, Vol. II, Prog. Math. {\bf 87} (1990), 111-195.
\bibitem
[Dr]{dr} V. Drinfeld, On Almost Cocommutative Hopf Algebras,
{\em Leningrad Mathematics Journal} {\bf 1} (1990), 321-342.
\bibitem
[DM]{dm} P. Deligne and J. Milne, Tannakian Categories, {\em Lecture
Notes in Mathematics} {\bf 900}, 101-228, 1982.
\bibitem
[EG1]{eg1} P. Etingof and S. Gelaki, Some Properties of
Finite-Dimensional Semisimple Hopf Algebras, {\em Mathematical
Research Letters} {\bf 5} (1998), 191-197.
\bibitem
[EG2]{eg2} P. Etingof and S. Gelaki, A Method of Construction of
Finite-Dimensional Triangular Semisimple Hopf Algebras, {\em
Mathematical Research Letters} {\bf 5} (1998), 551-561.
\bibitem
[EG3]{eg3} P. Etingof and S. Gelaki, The Representation Theory of
Cotriangular Semisimple Hopf Algebras, {\em International Mathematics
Research Notices} {\bf 7} (1999), 387-394.
\bibitem
[EG4]{eg4} P. Etingof and S. Gelaki, The Classification of Triangular
Semisimple and Cosemisimple Hopf Algebras Over an Algebraically Closed
Field, {\em International Mathematics Research Notices} {\bf 5} (2000),
223-234.
\bibitem
[EG5]{eg5} P. Etingof and S. Gelaki, On Finite-Dimensional Semisimple and
Cosemisimple Hopf Algebras In Positive Characteristic, {\em
International Mathematics Research Notices} {\bf 16} (1998), 851-864.
\bibitem
[EG6]{eg6} P. Etingof and S. Gelaki, On Cotriangular Hopf
Algebras, {\em American Journal of Mathematics}, to appear,
math.QA/0002128.
\bibitem
[EG7]{eg7} P. Etingof and S. Gelaki, Semisimple Hopf Algebras of
Dimension $pq$ are Trivial, {\em Journal of Algebra} {\bf 210}
(1998), 664-669.
\bibitem
[EG8]{eg8} P. Etingof and S. Gelaki, Classification of
Finite-Dimensional Triangular Hopf Algebras with the Chevalley
Property, {\em Mathematical Research Letters}, to appear,
math.QA/0101049.
\bibitem
[EG9]{eg9} P. Etingof and S. Gelaki, Isocategorical Groups, {\em
International Mathematics Research Notices} {\bf 2} (2001), 59-76.
\bibitem
[EGGS]{eggs} P. Etingof, S. Gelaki, R. Guralnick and J. Saxl,
Biperfect Hopf Algebras, {\em Journal of Algebra} {\bf 232} {\bf
no.1} (2000), 331-335.
\bibitem
[ES]{es} P. Etingof and O. Schiffmann, Lectures on Quantum
Groups, Lectures in Mathematical Physics, {\em
International Press, Boston, MA} (1998).
\bibitem
[G1]{g1} S. Gelaki, Some Properties and Examples of Pointed
Triangular Hopf Algebras, {\em Mathematical Research Letters} {\bf
6} (1999), 563-572; see corrected version at math.QA/9907106.
\bibitem
[G2]{g2} S. Gelaki, On the Classification of Finite-Dimensional
Triangular Hopf Algebras, {\em MSRI Publications volume dedicated
to Hopf Algebras}, to appear.
\bibitem
[HI]{hi} R.B. Howlett and I.M. Isaacs, On Groups of Central Type,
{\em Mathematische Zeitschrift} {\bf 179} (1982), 555-569.
\bibitem
[Kap]{kap} I. Kaplansky, Bialgebras, University of Chicago,
1975.
\bibitem [Kas]{kas} C. Kassel, Quantum Groups, Springer, New
York, 1995.
\bibitem [LR1]{lr1} R. Larson and D. Radford,
Semisimple
Cosemisimple Hopf Algebras, {\em American Journal of Mathematics}
{\bf 110} (1988), 187-195.
\bibitem
[LR2]{lr2} R.G. Larson and D.E. Radford,
Finite-Dimensional
Cosemisimple Hopf Algebras in Characteristic $0$ are Semisimple, {\em J.
Algebra} {\bf 117} (1988), 267-289.
\bibitem
[Mac1]{mac1} S. Mac Lane, Natural associativity and commutativity, {\em Rice Univ.
Studies} {\bf 49} (1963), 28-46.
\bibitem
[Mac2]{mac2} S. Mac Lane, Categories for the working mathematician, {\em Graduate
Text in Mathematics} {\bf 5}, Springer Verlag, NY, 1971.
\bibitem
[Mas]{ma} A. Masuoka, Self dual Hopf algebras of dimension $p^3$ obtained by
extension, {\em J. Algebra} {\bf 178} (1995), 791-806.
\bibitem
[Mon1]{mon1} S. Montgomery, Hopf algebras and their actions
on rings, {\em CBMS Lecture Notes} {\bf 82}, AMS, 1993.
\bibitem
[Mon2]{mon2} S. Montgomery, Classifying finite-dimensional
semisimple Hopf algebras, Trends in the representation theory of
finite-dimensional algebras (Seattle, WA 1997), 265-275,
{\em Contemp. Math.} {\bf 229}, AMS, Providence, RI, 19978.
\bibitem
[Mov]{mov} M. Movshev, Twisting in group algebras of finite
groups, {\em
Func. Anal. Appl.} {\bf 27} (1994), 240-244.
\bibitem
[R]{r} D.E. Radford, Minimal quasitriangular Hopf algebras,
{\em Journal of Algebra} {\bf 157} (1993), 285-315.
\bibitem
[S]{s} M. Sweedler, Hopf Algebras, Benjamin Press, 1968.
\bibitem
[Z]{z} Y. Zhu, Hopf algebras of prime dimension, {\em International
Mathematical Research Notices} {\bf 1} (1994), 53-59.
\end{thebibliography}
\end{document}